\documentclass[12pt]{article}
\usepackage{amssymb}
\usepackage{amsfonts}
\begin{document}

{\LARGE \bf How to Define Global Lie Group \\ \\ Actions on Functions} \\ \\

{\bf  Elem\'{e}r E ~Rosinger} \\

Department of Mathematics \\
and Applied Mathematics \\
University of Pretoria \\
0002 Pretoria \\
South Africa \\
eerosinger@hotmail.com \\ \\

{\bf Abstract} \\

A particularly easy, even if for long overlooked way is presented for defining globally
arbitrary Lie group actions on smooth functions on Euclidean domains. This way is based on the
appropriate use of the usual parametric representation of functions. As a further benefit of
this way, one can define large classes of genuine Lie semigroup actions. Here "genuine" means
that, unlike in the literature, such Lie semigroups need no longer be sub-semigroups of Lie
groups, and instead, can contain arbitrary noninvertible smooth functions on Euclidean
domains. \\ \\

{\bf 1. Introduction} \\

The advantages of being able to define {\it global} actions for arbitrary Lie groups are well
known for at least six decades by now, as presented systematically in the celebrated text of
Chevalley, for instance. \\
Yet, even in the case of Lie groups acting on Euclidean spaces, and not on manifolds in
general, the customary approach has not been able to go beyond a mere local definition, when
it comes to actions on functions by arbitrary Lie groups, see for instance Bluman \& Kumei,
Ibragimov, or Olver [1-3]. \\

Rather surprisingly, this failure to define globally the action on functions of arbitrary Lie
groups is due to an elementary difficulty, which can easily be overcome by a {\it parametric}
definition of functions, as shown for the first time in Rosinger [6], see also Rosinger [7]. \\

This parametric approach proves to have in fact two important advantages, namely, one of {\it
calculus}, and the other of {\it functorial} nature. The calculus advantage relates to the
simple and well known fact that the partial derivatives of any order of a parametrically given
function can be computed from it, without first having to bring the function to the usual,
nonparametric form. The functorial advantage, relating perhaps even to a simpler fact, is the
one which will actually allow the most easy, direct and natural global definition of arbitrary
Lie group actions on functions. In fact, as shown in Rosinger [6,7] and mentioned in the
sequel, it allows as well for the equally easy global definition of a far larger class of Lie
{\it semigroup} actions. \\

As a general remark about the parametric approach to the global definition of arbitrary Lie
group actions on functions, it is rather ironic to note that, in an embryonic, partial and
{\it local} manner, this approach has in fact been in use for a long time by now. \\

Indeed, suppose given a smooth function $f : \Omega \longrightarrow \mathbb{R}$, with
$\Omega \subseteq \mathbb{R}^n$ nonvoid, open. Further, suppose given an arbitrary Lie group
$G$ acting on $M = \Omega \times \mathbb{R}$ according to \\

$~~~~~~ G \times M \longrightarrow M $ \\

Then the usual way this Lie group action on $M$ is extended to such functions $f$, and thus to
${\cal C}^\infty ( \Omega, \mathbb{R} )$, is as follows. We consider the graph of $f$, that is,
the set \\

$~~~~~~ \gamma_f = \{~ (x,f(x)) ~|~ x \in \Omega ~\} \subseteq M $ \\

Therefore, for any $g \in G$, we can define point-wise the action $g \gamma_f$ and obtain
again a subset of $M$. \\

Unfortunately however, in general, it will {\it not} be true that \\

$~~~~~~ g \gamma_f = \gamma_h $ \\

for a certain smooth function $h : \Omega \longrightarrow \mathbb{R}$, which function $h$ if
it existed, it would obviously correspond to the {\it global} action of $g$ on $f$, that is,
we would have \\

$~~~~~~ g f = h $ \\

And then, the usual way to define arbitrary Lie group actions on functions overcomes this
difficulty at the cost of no less than a {\it double localization}, Olver [1-3], namely

\begin{itemize}

\item $g$ is {\it restricted} to a neighbourhood of the identity $e \in G$, and in addition

\item $f$ is {\it restricted} to suitable nonvoid, open subsets $\triangle$ of $\Omega$.

\end{itemize}

It is clear, however, that the consideration of the graph $\gamma_f$ of $f$ amounts to
replacing $f : \Omega \longrightarrow \mathbb{R}$ by the following special {\it parametric}
form of it, see (3.3), (3.4) in the sequel, namely $f_* : \Omega \longrightarrow M$, where
$\Omega \ni x \longmapsto f_*(x) = (x,f(x)) \in M$. Furthermore, in this case $g \gamma_f$ is
nothing else but $g f_*$, that is, the action of $g$ on $f_*$, which can {\it always} be
defined globally, irrespective of the function $f$, or of the Lie group action $G$ on $M$. \\

Thus it becomes clear that the {\it only} difficulty we have ever faced when trying to define
globally arbitrary Lie group actions on functions is {\it not} at all related to Lie groups or
functions, but solely to our rather unformulated, and yet quite implacable {\it intent} to
have $g f_*$ retranslated into a usual, {\it nonparametric} function $h : \Omega
\longrightarrow \mathbb{R}$. \\

On the other hand, the parametric approach to Lie group actions introduced in Rosinger [6], is
adopted and pursued in its {\it full} extent, that is, without any sort of localization, this
being the simple and fundamental reason for the fact that {\it arbitrary} Lie group actions
can be defined {\it globally} on smooth functions. \\
Furthermore, as shown in Rosinger [6], this possibility to define globally arbitrary Lie group
actions on smooth functions can easily be extended to actions on large classes of generalized
functions, and in particular, distributions, one of the effects of such an extension being the
first general solution of Hilbert's Fifth Problem, Rosinger [6]. \\
Also as mentioned and shown briefly in the sequel, one can define globally on functions the
action of far larger classes of Lie semigroups. This comes as a rather unexpected {\it bonus},
and the effect of the mentioned functorial nature of the parametric approach to Lie group
actions which allows the definition of arbitrary smooth - thus typically {\it noninvertible} -
actions. Such noninvertible actions can, of course, no longer belong to Lie group actions, but
only to Lie {\it semigroup} actions, Rosinger [6,7]. \\

Let us mention here in passing that the interest in such Lie semigroups of actions comes from
the fact that they range over a significantly {\it larger} class of actions than those
corresponding to Lie groups. Therefore, when applied to the study of solutions of PDEs - this
time as {\it semisymmetries} - they can offer new additional insights. \\
Furthermore, as pointed out by P J Olver, semigroups of actions appear quite naturally in
several aspects of the classical Lie theory, see for details Rosinger [6, chap. 13], [7]. \\ \\

{\bf 2. Difficulties with Actions on Usual Functions} \\

{\bf Classical Lie Group Actions}. For convenience, let us consider the familiar and important
setup when Lie group actions are used in the study of PDEs. In such cases, we are given a
linear or nonlinear PDEs of the general form \\

(2.1)~~~ $ T(x,D)U(x) ~=~ 0,~~~x \in \Omega ~\subseteq~ \mathbb{R}^n $ \\

where $\Omega$ is nonvoid open, $U :\Omega \longrightarrow \mathbb{R}$ is the unknown function,
while $T(x,D)$ is a ${\cal C}^\infty$-smooth linear or nonlinear partial differential operator.
The relevant Lie groups $G$ act on the open subset $M = \Omega \times \mathbb{R} \subseteq
\mathbb{R}^{n + 1}$, according to \\

(2.2)~~~ $ G \times M \ni (g, (x,u)) \longmapsto g(x,u) =
                                              (g_1(x,u), g_2(x,u)) \in M $ \\

where $x \in\Omega,~ u \in\mathbb{R}$ are the independent and dependent variables,
respectively, and \\

(2.3)~~~ $ \begin{array}{l}
                 G \times M \ni~(g,(x,u)) \longmapsto~ g_1(x,u) \in \Omega \\
                 G \times M \ni~(g,(x,u)) \longmapsto~ g_2(x,u) \in \mathbb{R}
            \end{array} $ \\

with $g_1$ and $g_2$ being ${\cal C}^\infty$-smooth. \\

We note that, given $g \in G$, in view of the Lie group axioms, it follows that the mapping \\

(2.4)~~~ $ M \ni (x,u) \stackrel{g}{\longmapsto} g(x,u) \in M $ \\

is a ${\cal C}^\infty$-smooth diffeomorphism. \\

A first basic problem in Lie group theory, when applied to PDEs, is how to {\it extend} the
action in (2.2), (2.3) of the Lie group $G$ on the open subset $M$, to an action of $G$ on the
${\cal C}^\infty$-smooth functions \\

(2.5)~~~ $ U : \Omega \longrightarrow \mathbb{R} $ \\

or more generally, on ${\cal C}^\infty$-smooth functions \\

(2.6)~~~ $U : \triangle \longrightarrow \mathbb{R} $ \\

where $\triangle \subseteq \Omega$ is nonvoid, open. And unless one solves this problem, one
simply cannot speak about the Lie group invariance of classical solutions of PDEs. \\

From this point of view, the Lie group actions (2.2), (2.3) are divided in two types, Olver
[1,2]. \\

The simpler ones, called {\it projectable}, or {\it fibre preserving}, satisfy the condition,
see (2.3) \\

(2.7)~~~ $ g_1(x,u) ~=~ g_1(x),~~~g \in G,~ (x,u) \in M $ \\

The special interest in Lie group actions (2.7) comes from the fact that they allow an easy
{\it global} extension to action on ${\cal C}^\infty$-smooth functions. Indeed, in this case,
in view of (2.4), it follows that for $g \in G$, we obtain the ${\cal C}^\infty$-smooth
diffeomorphism \\

(2.8)~~~ $ \Omega \ni x \stackrel{g_1}{\longmapsto} g_1(x) \in \Omega $ \\

Now, given $g \in G$ and $U$ in (2.6), it is easy to define the respective global Lie group
action \\

(2.9)~~~ $ g~U ~=~ \widetilde{U} : \widetilde{\triangle} ~=~ g_1(\triangle) \longrightarrow \mathbb{R} $ \\

by \\

(2.10)~~~ $ \widetilde{U} (g_1(x)) ~=~ g_2(x, U(x)),~~~x \in \triangle $ \\

Indeed, (2.4) implies that in (2.9), we have $\widetilde{\triangle} \subseteq \Omega$ nonvoid,
open, while (2.10) is equivalent with \\

(2.11)~~~ $ \widetilde{U}(\widetilde{x}) ~=~ g_2 (g_1^{-1}(\widetilde{x}),
                     U(g_1^{-1}(\widetilde{x}))),~~~\widetilde{x} \in \widetilde{\triangle}$ \\

However, an arbitrary Lie group action (2.2), (2.3) need {\it not} be projectable. And in such
a case the global extension of the Lie group action (2.2), (2.3) to ${\cal C}^\infty$-smooth
functions (2.5), or in general (2.6), will typically {\it fail}. In this way, we are obliged,
Olver [1,2], to limit ourselves to {\it local} Lie group actions on functions, and thus return
to the pre-Chevalley stage of Lie group theory. \\

Indeed, in the case of general, nonprojectable Lie group actions (2.2), (2.3), we may
immediately run into the problem of possible {\it noninvertibility}. Namely, certain ${\cal
C}^\infty$-smooth mappings involved in the definition of the group action $g~U = \widetilde{U}
: \widetilde{\triangle} \longrightarrow \mathbb{R}$ may fail to have inverses, let alone,
${\cal C}^\infty$-smooth ones. Let us illustrate this phenomenon in more detail. Given $g \in
G,$ let us write (2.3) in the form \\

(2.12)~~~ $ \begin{array}{l}
                 \widetilde{x} ~=~ g_1 (x,u) \\
                 \widetilde{u} ~=~ g_2 (x,u)
            \end{array} $ \\

where $(x,u),~ (\widetilde{x}, \widetilde{u}) \in M.$ Given now $U : \triangle \longrightarrow
\mathbb{R}$ as in (2.6), the natural way to define the group action $g U = \widetilde{U} :
\widetilde{\triangle} \longrightarrow \mathbb{R}$ would be by the relation, see (2.12) \\

(2.13)~~~ $ \widetilde{U} (g_1 (x,U(x))) ~=~ g_2 (x,U(x)),~~~x \in \triangle $ \\

which means that $\widetilde{U}(\widetilde{x}) = \widetilde{u}$.  However, in order that
(2.13) be a correct definition, we have to be able to obtain $x \in \triangle$ as a ${\cal
C}^\infty$-smooth function of $\widetilde{x} \in \widetilde{\triangle}$, by using the first
equation in \\

(2.14)~~~ $ \begin{array}{l}
                  \widetilde{x} ~=~ g_1(x,U(x)) \\
                  \widetilde{u} ~=~ g_2 (x, U(x))
            \end{array} $ \\

and thus by replacing $x \in \triangle$ in the second equation above, in order to obtain
$\widetilde{u}$ as a function of $\widetilde{x}$, that is, the relation (2.13). Furthermore,
one also has to obtain $\widetilde{\triangle} \subseteq~ \Omega$ as being nonvoid, open. The
crucial issue here is, therefore, the ${\cal C}^\infty$-{\it smooth invertibility} of the
mapping \\

(2.15)~~~ $ \triangle \ni x \stackrel{\alpha}{\longmapsto} g_1 (x,U(x)) \in \Omega $ \\

which obviously depends on $g$ and $U$. And as seen in the very simple example next, this in
general is not possible. \\

{\bf Example 2.1.} \\

Let us consider the following {\it nonprojectable} case of the Lie group action (2.2), (2.3),
where $\Omega = \mathbb{R},~ M = \Omega \times \mathbb{R} = \mathbb{R}^2,~ G = (\mathbb{R},
+ )$, and for $\epsilon = g \in G = \mathbb{R},~(x, u) \in M$, we have \\

$~~~ \begin{array}{l}
           \widetilde{x} ~=~ x + \epsilon u^2 \\
           \widetilde{u} ~=~ u
     \end{array} $ \\

Let us take $\triangle = \Omega = \mathbb{R}$ and the simple function $U : \triangle
\longrightarrow \mathbb{R}$ defined by $U(x) = x$, with $x \in \triangle$.  Then (2.15)
becomes \\

$~~~ \mathbb{R} \ni x \stackrel{\alpha}{\longmapsto} x + \epsilon x^2 \in \mathbb{R} $ \\

which is {\it not} invertible as a function, let alone as a ${\cal C}^\infty$-smooth function,
except for the trivial group action corresponding to $\epsilon = 0$, that is, to the identical
group transformation.

\hfill $\Box$ \\

The usual way to deal with this situation, Olver [1,2], is to consider the group action (2.2),
(2.3) as well as the mapping $\alpha$ in (2.15), and therefore the function to be acted upon
$U : \triangle \longrightarrow \mathbb{R}$, {\it only locally}, that is, to restrict all of
them to such suitable neighbourhoods of the neutral element $e \in G$, as well as of points
$x \in \triangle,$ on which $\alpha$ is ${\cal C}^\infty$-smooth invertible. \\

It is useful to note however that, depending also on the function $U$ in (2.6), the mapping
$\alpha$ in (2.15) can sometime happen to have a global, and not only local ${\cal
C}^\infty$-smooth inverse, even in the case of a nonprojectable Lie group action. For instance,
this happens if in the above Example 2.1., we consider $\widetilde{x} = x + \epsilon u$. \\

Let us mention what happens when the mapping $\alpha$ in (2.15) is invertible, regardless of
the Lie group action being projectable or not, and when its inverse $\alpha^{-1}$ is also a
${\cal C}^\infty$-smooth mapping. Then we can indeed turn to (2.13) in order to define the
group action $g~U = \widetilde{U}$ by \\

(2.16)~~~ $ \widetilde{U}~(\widetilde{x}) ~=~ g_2 (\alpha^{-1}(\widetilde{x}),
           U(\alpha^{-1}(\widetilde{x}))),~~~ \widetilde{x} \in~\widetilde{\triangle}$ \\

where \\

(2.17)~~~ $ \widetilde{\triangle} ~=~ \alpha (\triangle) ~~\mbox{is open} $ \\

Obviously, the case of {\it projectable} Lie group actions in (2.7) - (2.11) is included in
(2.16), (2.17). \\

As mentioned in the Introduction, here, following Rosinger [6,7], we take a {\it new route},
when dealing with the difficulties in (2.12) - (2.15), which we face in the case of general,
{\it nonprojectable} Lie group actions (2.2), (2.3). This new route will {\it not} require the
above mentioned traditional localisation of $g \in G,~ \alpha$ or $U$. In other words, we are
able to perform {\it globally} arbitrary Lie group actions on functions $U$ defined on the
whole of their unrestricted, original domains, as for instance in (2.5) and (2.6). Fortunately,
this construction is particularly simple and applicable without any undue restrictions. \\

{\bf A Simple, Basic Observation}. To summarize. The basis upon which we can delvelop this
global approach is the following rather simple observation :

\begin{itemize}

\item The usual impediment which prevents us from extending arbitrary Lie group actions (2.2),
(2.3) to global actions on functions (2.5) or (2.6) is {\it not} at all related to Lie groups,
but to the usual way of representing functions, by discriminating between independent and
dependent variables. Once one does away with such a discrimination, by using a parametric
representation of functions, the way to a natural and easy global Lie group action on
functions is open.

\end{itemize}

{\bf Parametrisation} in its essence amounts to the following {\it embedding} of the usual
definition of a function into a larger concept. Namely, a usual function \\

(2.18)~~~ $ A \ni x \stackrel{f}{\longmapsto} y ~=~ f(x) \in B $ \\

is actually {\it constrained} to be a correspondence from the set $A$ of its independent
variable $x$, to the set $B$ of its dependent variable $y$. \\

On the other hand, a parametric representation of $f$ can be given by any {\it pull-back} type
mapping \\

(2.19)~~~ $ P \ni p \stackrel{h}{\longmapsto} h(p) ~=~ (x(p),y(p)) \in A \times B $ \\

which maps any suitably given parameter domain $P$ into the graph of $f$, under the following
two conditions : \\

(2.20)~~~ $ y(p) ~=~ f(x(p)),~~~ p \in P $ \\

and \\

(2.21)~~~ $ P \ni p \longmapsto x(p) \in A ~~\mbox{is surjective} $ \\

With respect to $P$, this, in general, only implies that its cardinal is not
smaller than that of $A$. \\

However, when dealing with Lie group actions, the parameter domain $P$ is required to be a
suitable open subset in an Euclidean space, while the parametrisation $h$ is assumed ${\cal
C}^\infty$-smooth. \\

It follows that, in general, a parametric representation will introduce an {\it additional}
variable $p$, ranging over $P$, which this time is mapped into the pair $(x(p),y(p))$ of the
original independent and dependent variables, pair which is an element in the cartesian
product $A \times B$. \\

This kind of embedding, obtained by introducing an additional variable, and thus going beyond
the constraint of only dealing with the usual independent and dependent variables, proves to
have an important and naturally built in advantage. Namely, it allows for the first time - and
in a straightforward manner - the {\it global} definition of arbitrary Lie group actions on
functions. \\

In the usual, that is, nonparametric approach, however, when one wanst to define the Lie
group action on a function, and obtain again a function, one cannot in general do so, unless
at the end one is able to {\it separate} the independent and dependent variables, by
expressing the latter as a function of the former. And in the nonprojectable case of Lie group
actions, this typically is not possible, except locally in the independent variable, and also,
near to the trivial, identical Lie group transformation. \\

On the other hand, if one starts, and ends, with parametrically given functions, then as shown
in Rosinger [6,7] and seen in the sequel, one has no difficulties at all. \\ \\

{\bf 3. Parametric Functions} \\

{\bf Need for a Global Approach}. It is instructive to give another simple example, which by
its particularly familiar setup, can further highlight the rather basic, yet extreme
difficulties one may face when trying to define {\it globally} the action of a {\it
nonprojectable} Lie group on a function. \\

{\bf Example 3.1.} \\

Let us consider the Lie group action given by the usual {\it rotation of the plane}. In terms
of (2.2), (2.3), it means that $\Omega = \mathbb{R},~ M = \Omega \times \mathbb{R} =
\mathbb{R}^2,~ G = (\mathbb{R}, + )$ and for $\theta = g \in G = \mathbb{R},~ (x,u) \in M$, we
have \\

$~~~ \begin{array}{l}
            \widetilde{x} ~=~ x \cos \theta - u \sin \theta \\
            \widetilde{u} ~=~ x \sin \theta + u \cos \theta
      \end{array} $ \\

therefore, here again, we are dealing with a nonprojectable Lie group action, since for a
given $\theta = g$, obviously $\widetilde{x}$ depends not only on $x$, but also on $u$, see
(2.7). \\

Let $\triangle = \Omega = \mathbb{R},$ and $U : \triangle \longrightarrow \mathbb{R}$ be given
by the {\it parabola} $U(x) = x^2$, with $x \in \triangle$. \\

Then (2.15) takes the form \\

$~~~ \mathbb{R} \ni x \stackrel{\alpha}{\longmapsto}
                          x \cos \theta - x^2 \sin \theta \in \mathbb{R} $ \\

which, again, is {\it not} invertible, except for the trivial group actions, for which $\theta
= k\pi$, with $k \in \mathbb{Z}$. \\

It follows that, except for a trivial rotation of $\theta = \pm \pi$, which in this case
amounts to nothing else but a mere symmetry with respect to the origin of coordinates, the
parabola \\

$~~~ \triangle = \mathbb{R} \ni x \longmapsto x^2 \in \mathbb{R} $ \\

when taken as a whole, {\it cannot} be rotated at all in the plane, without ceasing to be the
{\it graph} of a function from $\triangle = \mathbb{R}$ to $\mathbb{R}$. Yet it is clear that,
as a geometric object, by arbitrarily rotating in the plane a parabola, one again gets a
parabola. \\

Therefore, the difficulty must lie with the particular way one happens to {\it represent} the
parabola, that is, as a function from $\triangle = \mathbb{R}$ to $\mathbb{R}$. \\

{\bf Parametric Representations}. It turns out that the alternative way to represent functions
$U : \triangle \longrightarrow \mathbb{R}$ in (2.5), (2.6), namely, {\it parametrically},
avoids the above difficulties related to the possible {\it lack of} ${\cal C}^\infty$-{\it
smooth inverse} of the mapping in (2.15), thus allows for the definition of global Lie group
actions on the resepctive functions. \\

Let us recall here that parametric representation, and not only of functions, is a rather
familiar method in differential geometry, among others, where it is used to define, for
instance, the concept of submanifold. \\
Here, parametric representation is only employed for functions such as those in (2.5),
(2.6). \\

Given therefore a ${\cal C}^\infty$-smooth function \\

(3.1)~~~ $ U : \triangle \longrightarrow \mathbb{R} $ \\

where $\triangle \subseteq \mathbb{R}^n$  is nonvoid, open, we denoted its {\it graph} by \\

(3.2)~~~ $ \gamma_{_U} ~=~ \{~(x, U(x)) ~|~ x \in \triangle ~\} ~\subseteq~ M $ \\

Now, a {\it parametric} representation of $U$ is given by {\it any} ${\cal C}^\infty$-smooth
function \\

(3.3)~~~ $ V : \Lambda ~\longrightarrow~ M $ \\

where the set $\Lambda \subseteq \mathbb{R}^n$ of {\it parameters} is nonvoid, open, and such
that \\

(3.4)~~~ $ V ( \Lambda ) ~=~ \gamma_{_U} $ \\

As seen in (3.22) - (3.25) and (3.31) - (3.34) in the sequel, the above condition (3.4),
although seemingly quite weak, has certain useful implications. \\

It is important to note that in (3.3), the set $\Lambda$ of parameters is assumed to be
n-dimensional. This however, is in line with the fact that the domain of definition
$\triangle$ of $U$ in (3.1) is also n-dimensional. In particular, since $M$ in (3.2) is
n+1-dimensional, it follows that condition (3.4) is quite natural. Later, when in (3.6) we
define the class of arbitrary parametrically given functions, which are of interest here, we
shall hold to this assumption on the dimension of the set of parameters. \\

{\bf Canonical Parametrisations}. Clearly, an immediate, simple and natural parametric
representation of $U$ in (3.1) is given by \\

(3.5)~~~ $ \triangle \ni x ~\longmapsto~ U_* ( x ) ~=~ ( x, U ( x ) ) \in M $ \\

and we shall call ${U_*} : \triangle \longrightarrow M$ the {\it canonical} parametric
representation of $U : \triangle \longrightarrow \mathbb{R}$. Thus in terms of (3.3), we have
$\Lambda = \triangle$ and $V = U_*$, and clearly, condition (3.4) is satisfied. \\

However, it is obvious that a function $U$ in (3.1) can have many other parametric
representations (3.3), (3.4).  Details in this respect can be found in the sequel. In
particular, we shall see in (3.22) and (3.24) that in a certain sense $U_*$ is the {\it
simplest possible} parametric representation of $U$. \\

{\bf Classes of Parametrisations}. Clearly, the set of functions in (3.3), (3.4) is {\it
larger} than that in (3.1). More precisely, not every function $V$ in (3.3), (3.4) is the
parametric representation of a function $U$ in (3.1). For instance, a nontrivially rotated
parabola can be written as a function in (3.3), (3.4), but not as a function in (3.1). \\

Let us, therefore, denote by \\

(3.6)~~~ $ {\cal C}_n^\infty(M) $ \\

the set of all ${\cal C}^\infty$-smooth functions $V : \Lambda \longrightarrow M$, where
$\Lambda \subseteq \mathbb{R}^n$ is nonvoid, open, and call them n-dimensional {\it parametric
representations} in $M$. Also, let us denote by \\

(3.7)~~~ $ {\cal C}_{par}^\infty(\Omega) $ \\

the set of ${\cal C}^\infty$-smooth {\it partial} functions $U : \triangle \longrightarrow
\mathbb{R}$, see (3.1), where $\triangle \subseteq \Omega$ is nonvoid, open. Then (3.5) yields
an {\it embedding} \\

(3.8)~~~ $ {\cal C}_{par}^\infty (\Omega) \ni U \longmapsto U_* \in {\cal C}_n^\infty (M) $ \\

while on the other hand, in view of (3.3), (3,4), we obtain a mapping \\

(3.9)~~~ $ {\cal C}_{par}^\infty (\Omega) \ni U \longmapsto
                           {\cal P}_U \subseteq {\cal C}_n^\infty (M) $ \\

where ${\cal P}_U$ is the set of mappings $V$ in (3.3), which satisfy (3.4). In other words,
${\cal P}_U$ is the set of all parametric representations of $U$. And in view of (3.5), it is
clear that \\

(3.10)~~~ $ U_* \in {\cal P}_U \neq \phi,~~~~ U \in~{\cal C}_{par}^\infty ( \Omega ) $ \\

The {\it important} point is that the construction of {\it arbitrary nonlinear} Lie group
actions on the set of functions ${\cal C}_n^\infty (M)$ will {\it no longer} suffer from the
above difficulties related to the possible {\it lack} of a ${\cal C}^\infty$-{\it smooth
inverse} of the mappings in (2.15). \\

Similar to (3.7), it will be useful, for $\ell \in {\bf N},~ \ell \geq 1$ and $N \subseteq
\mathbb{R}^\ell$ nonvoid, open, to denote by \\

(3.11)~~~ $ {\cal C}_{par}^\infty (\Omega, N) $ \\

the set of all partial functions $U : \triangle \longrightarrow N$ which are ${\cal
C}^\infty$-smooth, where $\triangle \subseteq \Omega$ is any nonvoid, open subset. \\

Obviously ${\cal C}^\infty_{par}(\Omega, M) \subseteq {\cal C}^\infty_n(M)$. \\

{\bf Comparing Parametrisations}. Here we further clarify the meaning of the parametric
representation of functions, defined in (3.1) - (3.6). \\

Let us define a {\it preorder} $\leq$ on ${\cal C}_n^\infty (M),$ that is, a reflexive and
transitive binary relation, as follows. Given $\Lambda \stackrel{V}{\longrightarrow} M$ and
$\Lambda^ \prime \stackrel{V^\prime}{\longrightarrow} M$, with $\Lambda, \Lambda^ \prime
\subseteq \mathbb{R}^n$ nonvoid, open, then \\

(3.12)~~~ $ V \leq V^ \prime $ \\

if and only if there exists a surjective ${\cal C}^\infty$-smooth mapping $\Lambda
\stackrel{\varphi}{\longrightarrow} \Lambda^ \prime$, such that the diagram is commutative

\bigskip
\begin{math}
\setlength{\unitlength}{1cm}
\thicklines
\begin{picture}(8,6)

\put(0,2.5){(3.13)}
\put(3,2.5){$\varphi$}
\put(3.65,5){\vector(0,-2){5}}
\put(3.5,5.6){$\Lambda$}
\put(4.5,5){\vector(3,-2){3}}
\put(6.1,4.1){$V$}
\put(8,2.5){$M$}
\put(4.5,0){\vector(3,2){3}}
\put(6.1,0.5){$V^\prime$}
\put(3.5,-0.7){$\Lambda^\prime$}
\end{picture}
\end{math}

\vspace{1.5cm}
It is easy to see that, owing to the surjectivity of $\varphi$, we obtain \\

(3.14)~~~ $ V(\Lambda) ~=~ V^ \prime(\Lambda^ \prime) $ \\

A natural interpretation of this preorder $V \leq V^ \prime$ is that the parametrisation
$V^ \prime$ is {\it simpler} that $V$. This is illustrated in \\

{\bf Example 3.2.} \\

Let $\Lambda = \Omega \subseteq \mathbb{R}^n$ be nonvoid, open, and let us take any
$\varphi : \Omega \longrightarrow \Omega$ which is ${\cal C}^\infty$-smooth and surjective,
but it is not injective. Also, let us take any ${\cal C}^\infty$-smooth $U : \Omega
\longrightarrow \mathbb{R}$. We can now define the parametric function in ${\cal
C}^\infty_n(M)$, namely \\

(3.15)~~~ $ \Lambda ~=~ \Omega \ni~x \stackrel{V}{\longmapsto}
                            V(x) ~=~ (\varphi(x), U(\varphi(x))) \in M $ \\

Then clearly \\

(3.16)~~~ $ V \leq U_* $ \\

However \\

(3.17)~~~ $ U_* \not\le V $ \\

indeed, assume $\Omega \stackrel{\Psi}{\longrightarrow} \Lambda$ is surjective and ${\cal
C}^\infty$-smooth, and such that $U_* = V~\circ~\Psi$. Then $x = \varphi (\Psi(x)), x \in
\Omega$, which means that, contrary to the assumption, $\varphi$ is injective, since $\Psi$ is
surjective.

\hfill $\Box$ \\

Recalling that, see (3.5) \\

(3.18)~~~ $ \Omega \ni x \stackrel{U_*}{\longmapsto} (x, U(x)) \in M $ \\

and comparing it with (3.15), where $\varphi$ can be arbitrary under the mentioned assumptions,
it follows that a natural meaning of (3.16), (3.17) is that the canonical parametric
representation $U_*$ of $U$ is {\it simpler} than {\it all} the other parametric
representations of $U$, given by $V$ in (3.15), see also (3.22) and (3.24) below. \\

{\bf Basic Properties of Parametric Representations}. It is useful to note that the simple
looking condition (3.4) is precisely the one which leads to the situation in Example 3.2.
Indeed, let $U~:~\triangle \longrightarrow \mathbb{R}$ be a ${\cal C}^\infty$-smooth function,
with $\triangle \subseteq \Omega$ nonvoid, open, and let $\Lambda \stackrel{V}{\longrightarrow}
M$ be any function in ${\cal C}_n^\infty (M),$ which therefore acts according to \\

(3.19)~~~ $ \begin{array}{l}
                \Lambda \ni y \longmapsto  V(y) ~=~ (V_1(y), V_2(y)) \in M \\ \\
                \Lambda \ni y \longmapsto V_1(y) \in \Omega \\ \\
                \Lambda \ni y \longmapsto V_2(y) \in \mathbb{R}
            \end{array} $ \\

Then it is easy to see that the following four conditions are equivalent: \\

(3.20)~~~ $V$ is a parametric representation of $U$ \\

(3.21)~~~ $ V (\Lambda) ~=~ \gamma_U $ \\

(3.22)~~~ $ V ~\leq~ U_* $ \\

(3.23)~ $V_1$ is surjective and $V_2 = U \circ V_1$ \\

In view of the above and (3.9), it follows that \\

(3.24)~~~ $ \begin{array}{l}
                  \forall~~~ V \in {\cal P}_U ~: \\ \\
                  ~~~~~ V \leq U_*
            \end{array} $ \\

In this way, in view of (3.22) and (3.24), and in the sense of Example 3.2., it is clear that
for any given function $U$, its {\it canonical} parametric representation ${U_*}$ is {\it
simpler} than any other parametric representation of that function. \\

Moreover, given two ${\cal C}^\infty$-smooth functions $U_1, U_2 : \triangle \longrightarrow
\mathbb{R}$ in ${\cal C}_{par}^\infty (\Omega),$ then \\

(3.25)~~~ $ ( U_1 )_* ~\leq~ ( U_2 )_* ~~\Longrightarrow~~ U_1 ~=~ U_2 $ \\

which shows to what a large extent the canonical parametric representation does in fact
determine a function. \\

{\bf Staying with Usual Functions}. We conclude that given $\Lambda
\stackrel{V}{\longrightarrow} M$ in ${\cal C}_n^\infty (M),$ then a sufficient condition for
the existence of a function $U : \triangle \longrightarrow \mathbb{R}$ in ${\cal
C}_{par}^\infty (\Omega)$, such that $V$ is a parametric representation of $U,$ is given by,
see (3.23) \\

(3.26)~~~ $\Lambda \stackrel{V_1}{\longrightarrow} \triangle$ is a ${\cal C}^\infty$-smooth
                                                     diffeomorphism \\

In this case it also follows that, see (3.19) \\

(3.27)~~~ $ U ~=~ V_2 \circ V_1^{-1} $ \\

as well as, see (3.22) \\

(3.28)~~~ $ V \leq U_* $ \\

However, when one deala with arbitrary nonlinear, and possibly nonprojectable Lie group
actions on functions, one can encounter the general situation of mappings $\Lambda
\stackrel{V}{\longrightarrow} M$ in ${\cal C}_n^\infty (M)$ which may fail to satisfy (3.26).
Thus this condition (3.26) can be seen as the more general reformulation of the ${\cal
C}^\infty$-smooth invertibility problem in (2.15). \\

{\bf Equivalent Parametrisations}. Let us consider two parametric functions $V : \Lambda
\longrightarrow M$ and $V\,' : \Lambda\,' \longrightarrow M$ in ${\cal C}^\infty_n(M)$, see
(3.6). We say that $V$ and $V\,'$ are {\it equivalent}, and write \\

(3.29)~~~ $ V \approx V\,' $ \\

if and only if, see (3.4) \\

(3.30)~~~ $ V(\Lambda) ~=~ V\,'(\Lambda\,') $ \\

Clearly, if we have a usual function $U$ in (3.1), then in view of (3.4), $V$ will be a
parametric representation of $U$, if and only if, see (3.5) \\

(3.31)~~~ $ V \approx U_* $ \\

Also, owing to (3.12) - (3.14) and (3.30), it follows that \\

(3.32)~~~ $ V \leq V\,' ~~\Longrightarrow~~ V \approx V\,' $ \\ \\

{\bf 4. Actions on Parametric Representations} \\

{\bf Natural Definition}. Now with the use of parametric representations, we come to the {\it
basic idea} in this paper, namely, we can define the {\it arbitrary Lie group actions on
functions} \\

(4.1)~~~ $ G \times {\cal C}_n^\infty (M) \longrightarrow {\cal C}_n^\infty (M) $ \\

in the following simple and natural way. Given $g \in G$ and a function $\Lambda
\stackrel{V}{\longrightarrow} M$ from ${\cal C}_n^\infty (M)$, we define \\

(4.2)~~~ $ g~V ~=~ g \circ V $ \\

where in the right hand term, $g$ is the mapping in (2.4). In other words, we use as
definition of the Lie group action the commutative diagram

\bigskip
\begin{math}
\setlength{\unitlength}{1cm}
\thicklines
\begin{picture}(10,3)

\put(1.5,2){$\Lambda$}
\put(3.7,2.3){$V$}
\put(2,2.1){\vector(1,0){4}}
\put(6.3,2){$M$}
\put(8.7,2.35){$g$}
\put(7,2.1){\vector(1,0){4}}
\put(11.4,2){$M$}
\put(1.6,0.95){\line(0,1){0.75}}
\put(0,1.65){(4.3)}
\put(1.59,0.95){\line(1,0){10.05}}
\put(11.6,0.95){\vector(0,1){0.75}}
\put(6.3,0.4){$gV$}
\end{picture}
\end{math} \\

Clearly, with the definition (4.2), (4.3), we have \\

(4.4)~~~ $ g V \in {\cal C}_n^\infty (M),~~~ \Lambda \stackrel{g V}{\longrightarrow} M $ \\

that is, $g V$ and $V$ have the {\it same} domain of definition $\Lambda$, and the {\it same}
range $M$. \\

{\bf Properties}. We show now that the Lie group actions (4.1) contain as a particular case
the usual Lie group actions on functions, Olver [1,2], namely \\

$~~~~~~ G \times {\cal C}_{par}^\infty (\Omega) \longrightarrow
                                         {\cal C}_{par}^\infty (\Omega) $ \\

whenever the latter can be defined {\it globally}. Indeed, assume given $g \in G$ and $U :
\triangle \longrightarrow \mathbb{R}$ in (2.26), such that the mapping $\alpha$ in (2.15) is a
${\cal C}^\infty$-smooth diffeomorphism. In view of (3.5), we obtain $U_* \in {\cal
C}_n^\infty (M)$ and then (2.4) and (4.2) give \\

\begin{math}
\setlength{\unitlength}{1cm}
\thicklines
\begin{picture}(15,0.5)

\put(2.6,0.35){$gU_*$}
\put(0,0){$(4.5) ~~~~ \triangle$}
\put(2.15,0.15){\vector(1,0){1.5}}
\put(3.9,0){$M$}

\end{picture}
\end{math} \\

where \\

(4.6)~~~ $ (g~ U_*)(x) ~=~ g(U_*(x)) ~=~ g(x,U(x)) ~=~ (g_1(x,U(x)), g_2(x,U(x))) $ \\

with $x \in \triangle$.  On the other hand, in view of our assumption on $\alpha$, we can
apply (2.16), (2.17) and obtain \\

(4.7)~~~ $ (g~ U)(\widetilde{x}) ~=~ g_2 (\alpha^{-1}(\widetilde{x}),
              U(\alpha^{-1}(\widetilde{x}))),~~~ \widetilde{x} \in \widetilde{\triangle} ~=~
                                        \alpha(\triangle) $ \\

therefore (3.5) gives \\ \\

\begin{math}
\setlength{\unitlength}{1cm}
\thicklines
\begin{picture}(15,1.5)

\put(2.4,1.35){$(gU)_*$}
\put(1.5,1){$\widetilde{\triangle}$}
\put(0,0.5){(4.8)}
\put(2.2,1.15){\vector(1,0){1.5}}
\put(4.1,1){$M$}
\put(1.5,0){$\widetilde{x}$}
\put(2.2,0.15){\vector(1,0){1}}
\put(3.6,0){$(\widetilde{x},g_2(\alpha^{-1}(\widetilde{x}),U(\alpha^{-1}
             (\widetilde{x}))))$}

\end{picture}
\end{math} \\

Now, from (2.17), (4.5) and (4.8) it is clear that, in general \\

(4.9)~~~ $ g~ U_* ~\neq~ (g~ U)_* $ \\

since, for instance, their domains of definition need not be the same. However, we have \\

(4.10)~~~ $ (g~ U_*) (\triangle) ~=~ \gamma_{_gU} $ \\

since a direct computation gives, see (3.2), (2.2) \\

(4.11)~~~ $ \gamma_{_gU} ~=~ g \gamma_{_U} ~=~
                     \{~ (g_1(x,U(x)), g_2(x, U(x))) ~~|~~ x \in \triangle ~\} $ \\

and on the other hand, see (4.6), (2.2) \\

(4.12)~~~ $ (g~ U_*)(\triangle) ~=~
                     \{~ (g_1 (x, U(x)), g_2(x, U(x))) ~~|~~ x \in \triangle ~\} $ \\

It follows that, in view of (3.3), (3.4), the Lie group action $g~ U_*$ of $g$ on the
parametric representation $U_*$ of $U$, is itself a parametric representation of $g~ U$, which
is the Lie group action of $g$ on $U$.  In other words, in general, the diagram \\

\begin{math}
\setlength{\unitlength}{0.2cm}
\thicklines
\begin{picture}(60,22)

\put(10,19){${\cal C}_{par}^{\infty} (\Omega)  \ni~ U$}
\put(33,20.5){$g$}
\put(25,19){\vector(1,0){17}}
\put(45,19){$gU \in {\cal C}_{par}^{\infty} (\Omega)$}
\put(0,11){$(4.13)$}
\put(15,11){$(~~)_*$}
\put(20,17){\vector(0,-1){12}}
\put(46,17){\vector(0,-1){12}}
\put(48,11){$(~~)_*$}
\put(10,2){${\cal C}_n^{\infty} (M) \ni U_*$}
\put(25,2){\vector(1,0){10}}
\put(29,0){$g$}
\put(38,2){$gU_* \neq (g U)_* \in {\cal C}_n^{\infty} (M)$}

\end{picture}
\end{math} \\

is {\it not} commutative, see (4.9).  Nevertheless, we have, see (3.9), (3.24), (4.11),
(4.12) \\

(4.14)~~~ $ g~ U_* \in {\cal P}_{gU},~~~ g~U_* ~\leq~ (g~U)_* $ \\

Further, we note that, regardless of (4.9) and (4.13), we obtain the following {\it
commutative} diagram \\

\begin{math}
\setlength{\unitlength}{1cm}
\thicklines
\begin{picture}(8,6.5)

\put(0,3.2){(4.15)}
\put(2.8,3.2){$\alpha$}
\put(3.45,5.7){\vector(0,-2){5}}
\put(3.85,0.7){\vector(0,2){5}}
\put(4.3,3.2){$\alpha^{-1}$}
\put(3.5,6.3){$\triangle$}
\put(4.5,5.7){\vector(3,-2){3}}
\put(6.1,4.9){$gU_*$}
\put(8,3.2){$M$}
\put(4.5,0.7){\vector(3,2){3}}
\put(6.1,1.2){$(gU)_*$}
\put(3.5,0){$\triangle^\prime$}
\end{picture}
\end{math} \\

which follows easily from (4.6), (4.8) and (3.5).  In this way, in view of (3.12), (3.13), the
above commutative diagram means that \\

(4.16)~~~ $ g~U_* ~\leq~ (g~U)_* ~\leq~ g~U_* $ \\

in other words, in case the usual Lie group action $g~U$ of $g \in G$ on the function $U$
exists {\it globally}, then its {\it canonical} parametric representation $(g~U)_*$ is both
{\it more simple} and {\it less simple} than $g~U_*$, which is the Lie group action on the
{\it canonical} parametric representation $U_*$ of $U$, and which {\it always} exists. \\

{\bf Remark 4.1.} \\

In view of the commutative diagram (4.15), and the double inequality in (4.16), it appears to
be natural to use the {\it global} Lie group action $g~U_*$, which always exists, instead of
the Lie group action $g~U$, since as seen in section 2, the latter need not always exist. \\

In fact, the {\it essential interest} in using parametric representations is that we can {\it
abandon} $(g~U)_*$ in (4.9), (4.13) - (4.16), and instead, use $g~U_*$, which always exist
{\it globally}, and which also happens to be a parametric representation of $g~U$, whenever
the latter exists globally in the classical sense, Olver [1,2].

\hfill $\Box$ \\

Finally, related to the commutative diagram (4.15), and the double inequality (4.16), we have
the following additional {\it universality} type properties. Given $ \Lambda
\stackrel{V}{\longrightarrow} M$ a function from ${\cal C}_n^\infty (M)$, such that the
diagram of ${\cal C}^\infty$-smooth mappings \\

\begin{math}
\setlength{\unitlength}{1cm}
\begin{picture}(8,6.5)
\thicklines

\put(0,3.2){(4.17)}
\put(3,3.2){$\lambda$}
\put(3.65,5.7){\vector(0,-2){5}}
\put(3.5,6.3){$\triangle^\prime$}
\put(4.5,5.7){\vector(3,-2){3}}
\put(6.05,4.85){$(gU)_*$}
\put(8,3.2){$M$}
\put(4.5,0.7){\vector(3,2){3}}
\put(6.1,1.2){$gV$}
\put(3.5,0){$\Lambda$}
\end{picture}
\end{math} \\

is commutative, then \\

(4.18)~~~ $ V \circ \lambda \circ \alpha ~=~ U_* $ \\

Indeed, (4.15) - (4.17) yield \\

$~~~~~~ g~U_* ~=~ g~V \circ \lambda \circ \alpha $ \\

hence (4.18) follows from (4.2). Similarly, if the diagram of ${\cal C}^\infty$-smooth
mappings \\

\begin{math}
\setlength{\unitlength}{1cm}
\thicklines
\begin{picture}(8,6.5)

\put(0,3.2){(4.19)}
\put(3,3.2){$\lambda$}
\put(3.65,5.7){\vector(0,-2){5}}
\put(3.5,6.3){$\Lambda$}
\put(4.5,5.7){\vector(3,-2){3}}
\put(6.05,4.85){$gV$}
\put(8,3.2){$M$}
\put(4.5,0.7){\vector(3,2){3}}
\put(6.1,1.2){$(gU)_*$}
\put(3.5,0){$\triangle^\prime$}
\end{picture}
\end{math} \\

is commutative, then \\

(4.20)~~~ $ V ~=~ U_* \circ \alpha^{-1} \circ \lambda $ \\

Also, if the diagram of ${\cal C}^\infty$-smooth mappings

\begin{math}
\setlength{\unitlength}{1cm}
\thicklines
\begin{picture}(8,6.5)

\put(0,2.5){(4.21)}
\put(3,2.5){$\lambda$}
\put(3.65,5){\vector(0,-2){5}}
\put(3.5,5.6){$\triangle$}
\put(4.5,5){\vector(3,-2){3}}
\put(6.05,4.15){$gU_*$}
\put(8,2.5){$M$}
\put(4.5,0){\vector(3,2){3}}
\put(6.1,0.5){$gV$}
\put(3.5,-0.7){$\Lambda$}
\end{picture}
\end{math} \\ \\

is commutative, then \\

(4.22)~~~ $ V \circ \lambda ~=~ U_* $ \\

Finally, if the diagram of ${\cal C}^\infty$-smooth mappings \\

\begin{math}
\setlength{\unitlength}{1cm}
\thicklines
\begin{picture}(8,6.5)

\put(0,3.2){(4.23)}
\put(3,3.2){$\lambda$}
\put(3.65,5.7){\vector(0,-2){5}}
\put(3.5,6.3){$\Lambda$}
\put(4.5,5.7){\vector(3,-2){3}}
\put(6.05,4.85){$gV$}
\put(8,3.2){$M$}
\put(4.5,0.7){\vector(3,2){3}}
\put(6.1,1.2){$gU_*$}
\put(3.5,0){$\triangle$}
\end{picture}
\end{math} \\

is commutative, then \\

(4.24)~~~ $ V ~=~ U_* \circ \lambda $ \\

The commutative diagram (4.15), the double inequality (4.16) and the universal properties
(4.17) - (4.24) give both the explanation and remedy for the failures in (4.9) and
(4.13). \\ \\

{\bf 5. Comments} \\

The {\it novelty} of the extension of Lie group actions to parametric functions, as defined in
(4.2), (4.3), when compared with the usual one in section 2, becomes now clear. Indeed, in the
usual approach, Bluman \& Kumei, Ibragimov, Olver [1-3], one proceeds as follows. \\

First, one defines the Lie group action (2.2) on the set $M$ of independent and dependent
variables, respectively, $x \in \Omega$ and $u \in \mathbb{R}$. \\

Then as a second step, one extends this initial Lie group action to functions $U : \triangle
\longrightarrow \mathbb{R}$ in (2.6). \\
This extension is done by replacing $U$ with its graph $\gamma_{_U} \subseteq M$, and then
letting the Lie group act pointwise on $\gamma_{_U}$, seen as a subset of $M$. Certainly, for
every $g \in G$, we obtain in this way a well defined subset $g \gamma_{_U} \subseteq M$.
However, for {\it nonprojectable} Lie groups it need {\it not} in general happen that \\

(4.25)~~~ $ \begin{array}{l}
               \exists~~~ \widetilde{U}\in {\cal C}_{par}^{\infty} (\Omega) ~: \\ \\
               ~~~~~ g \gamma_U ~=~ \gamma_{\widetilde{U}}
             \end{array} $ \\

that is, the subset $g \gamma_U$ need {\it not} be the graph of {\it any} function
$\widetilde{U} : \widetilde{\triangle} \longrightarrow \mathbb{R}$, where
$\widetilde{\triangle} \subseteq \Omega$ is nonvoid, open. In this way, the usual method of
extending the Lie group action (2.2) from the set $M$ to functions $U$ in ${\cal
C}_{par}^\infty (\Omega)$, by using the graph $\gamma_U \subseteq M$ of $U$, is severely
limited in its globality, in the case of nonprojectable Lie groups. \\

The way out of this nonglobality impasse, as presented in this paper, is based on the
observation that the functions $U$ in ${\cal C}_{par}^\infty (\Omega)$ need {\it not} be seen
as being defined in terms which are necessarily {\it internal} or {\it confined} to the set
$M$ of {\it independent and dependent variables} $x \in \Omega$ and $u \in \mathbb{R}$,
respectively. \\

Indeed, by introducing the {\it parametric representation} of functions in ${\cal
C}_{par}^\infty (\Omega)$, as done in section 3, we can embed ${\cal C}_{par}^\infty (\Omega)$
into the space of parametric functions ${\cal C}_n^\infty (M)$, see (3.8). These parametric
functions have {\it arbitrary} domains, which are nonvoid, open in $\mathbb{R}^n$, however,
their {\it range} is {\it always} in the set $M$ of independent and dependent variables. And
as seen in (4.1) - (4.4), extending {\it globally} arbitrary Lie group actions (2.2) to
functions in ${\cal C}_n^\infty (M)$ is a rather simple and straightforward procedure, being
merely the {\it composition} of two mappings, each of which always exists. Furthermore, as
seen in (4.5) - (4.24), this extension contains as a particular case the usual way Lie group
actions (2.2) are extended to functions in ${\cal C}_{par}^\infty (\Omega, \mathbb{R})$,
whenever these latter extensions happen to exist globally. \\ \\

\newpage

{\bf 6. Semigroup Actions} \\

It is obvious that (4.3) remains valid, that is, the composition $g~V$ will still be a ${\cal
C}^\infty$-smooth function, and thus an element of ${\cal C}^\infty_n(M)$, even if the
mapping \\

(6.1)~~~ $ M \stackrel{g}{\longrightarrow} M $ \\

is no longer restricted to being given by the Lie group action (2.2), through the ${\cal
C}^\infty$-smooth diffeomorphism (2.4), but instead, it is simply an arbitrary ${\cal
C}^\infty$-smooth function. \\

In other words, (4.3) actually defines the following extension of (4.1) \\

(6.2)~~~ $ {\cal C}^\infty(M,M) \times {\cal C}^\infty_n(M)
                     ~~\longrightarrow~~ {\cal C}^\infty_n(M) $ \\

And since ${\cal C}^\infty(M,M)$ is a noncommutative {\it semigroup} with identity, and not a
group, when considered with the usual composition of functions, it is clear that (6.2) is a
vast extension of (4.1), no matter which would be the Lie group $G$ considered in (2.2). \\

Here, it should be noted that there has been an interest in certain Lie {\it semigroup} type
actions, Hilgert \& Neeb, Weinstein. None of them, however, aims anywhere near to the
generality of (6.2). Indeed, in Hilgert \& Neeb, which follows the work of the school of K H
Hofmann, the Lie semigroups considered must be {\it subsemigroups} of Lie groups, thus they
cannot include the semigroup ${\cal C}^\infty(M,M)$ which acts in (6.2). As for the concept of
grupoid, presented in the survey of Weinstein, it is similarly not capable of including the
mentioned semigroups which act in (6.2). \\

Needless to say, the extension of the {\it symmetry} concept from the framework of the Lie
group actions in (4.1), to that of the {\it vastly more general semigroup} actions in (6.2),
can be of a significant interest, among others, in the study of PDEs, even in the case of
their classical solutions. \\

A start in this direction was presented in Rosinger [6, chap. 13] and Rosinger [7].

\end{document}